\documentclass[times,sort&compress,3p]{elsarticle} 
\journal{Journal of Multivariate Analysis}
\usepackage{amsmath, amssymb, enumerate, amsthm}

\usepackage{color, multirow, booktabs}
\usepackage{amssymb}
\usepackage{amsmath}

\newtheorem{thm}{Theorem}
\newtheorem{lem}{Lemma}

\newtheorem{prp}{Proposition}

\newtheorem{remark}{Remark}

\def\al{{\alpha}}
\def\be{{\beta}}
\def\ga{{\gamma}}

\def\la{{\lambda}}
\def\si{{\sigma}}

\def\bbe{{\text{\boldmath $\beta$}}}

\def\bde{{\text{\boldmath $\delta$}}}
\def\bep{{\text{\boldmath $\varepsilon$}}}

\def\bom{{\text{\boldmath $\omega$}}}

\def\bmu{{\text{\boldmath $\mu$}}}
\def\bnu{{\text{\boldmath $\nu$}}}

\def\bnabla{{\text{\boldmath $\nabla$}}}

\def\sih{{\hat \si}}

\def\bbeh{{\widehat \bbe}}

\def\bmuh{{\widehat \bmu}}

\def\De{{\Delta}}
\def\Si{{\Sigma}}
\def\Ga{{\Gamma}}

\def\bDe{{\text{\boldmath $\De$}}}
\def\bSi{{\text{\boldmath $\Si$}}}

\def\a{{\text{\boldmath $a$}}}

\def\h{{\text{\boldmath $h$}}}

\def\j{{\text{\boldmath $j$}}}

\def\x{{\text{\boldmath $x$}}}
\def\y{{\text{\boldmath $y$}}}

\def\A{{\text{\boldmath $A$}}}
\def\B{{\text{\boldmath $B$}}}

\def\I{{\text{\boldmath $I$}}}
\def\J{{\text{\boldmath $J$}}}

\def\Q{{\text{\boldmath $Q$}}}

\def\V{{\text{\boldmath $V$}}}

\def\X{{\text{\boldmath $X$}}}
\def\Y{{\text{\boldmath $Y$}}}
\def\Z{{\text{\boldmath $Z$}}}

\def\Xb{{\overline \X}}

\def\Nc{{\cal N}}

\def\Re{{\mathbb{R}}}
\def\tr{{\rm tr\,}}

\def\vec{{\bf vec\,}}

\def\Ch{{\rm Ch}}

\def\[{{\text{\boldmath $[$}}}
\def\]{{\text{\boldmath $]$}}}

\def\zero{{\bf\text{\boldmath $0$}}}

\def\|{{\,|\,}}

\def\/{{\Bigr/\!\!}}

\def\1r{{\rm (1)}}
\def\2r{{\rm (2)}}
\def\3r{{\rm (3)}}
\def\4r{{\rm (4)}}
\def\5r{{\rm (5)}}

\def\non{{\nonumber}}

\def\bnuh{{\widehat \bnu}}

\def\bMu{{\boldmath{\text{$\mathcal{M}$}}}}
\def\bMuh{{\widehat \bMu}}
\DeclareMathOperator{\blockdiag}{block\ diag}

\begin{document}

\begin{frontmatter}

\title{Bayesian simultaneous estimation for means in $k$-sample problems}

\author[label1]{Ryo Imai}
\author[label2]{Tatsuya Kubokawa}
\author[label3]{Malay Ghosh}
\address[label1]{Graduate School of Economics, University of Tokyo, 
7-3-1 Hongo, Bunkyo-ku, Tokyo 113-0033, Japan}
\address[label2]{Faculty of Economics, University of Tokyo, 
7-3-1 Hongo, Bunkyo-ku, Tokyo 113-0033, Japan}
\address[label3]{Department of Statistics, University of Florida, 102 Griffin-Floyd Hall, Gainesville, FL 32611, USA}

\begin{abstract}
This paper is concerned with the simultaneous estimation of $k$ population means when one suspects that the $k$ means are nearly equal.
As an alternative to the preliminary test estimator based on the test statistics for testing hypothesis of equal means, we derive Bayesian and minimax estimators which shrink individual sample means toward a pooled mean estimator given under the hypothesis.
It is shown that both the preliminary test estimator and the Bayesian minimax shrinkage estimators are further improved by shrinking the pooled mean estimator.
The performance of the proposed shrinkage estimators is investigated by simulation. 
\end{abstract}


\end{frontmatter}

\section{Introduction}

Consider the multivariate $k$-sample problem expressed in the following canonical form:
$p$-variate random vectors $\X_1, \ldots, \X_k $ and a positive scalar random variable $S$ 
are mutually independent and distributed as
\begin{equation}
\X_1\sim \Nc_p(\bmu_1, \si^2\V_1), \quad \ldots, \quad \X_k\sim \Nc_p(\bmu_k, \si^2\V_k), \quad
S/\si^2\sim \chi_{(n)}^2,
\label{eqn:model}
\end{equation}
where the $p \times 1$ means $\bmu_1, \ldots, \bmu_k$ and the scale parameter $\si^2$ are unknown, and $\V_1, \ldots, \V_k$ are $p\times p$ known and positive definite symmetric matrices.
In this model, we want to estimate $\bmu_1, \ldots, \bmu_k$ simultaneously relative to the quadratic loss function
\begin{equation}
L(\bde, \bom) = \frac{1}{\si^2} \sum_{i=1}^k \Vert\bde_i-\bmu_i\Vert_{\Q_i}^2 
= \frac{1}{\si^2} \sum_{i=1}^k (\bde_i-\bmu_i)^\top\Q_i(\bde_i-\bmu_i) ,
\label{eqn:loss}
\end{equation}
where $\Q_1, \ldots, \Q_k$ are $p\times p$ known and positive definite symmetric matrices, $\Vert \a \Vert_\A^2=\a^\top \A \a$ with $\a^\top$ standing for the transpose of $\a$, $\bom=(\bmu_1, \ldots, \bmu_k, \si^2)$ is a vector of unknown parameters and $\bde=(\bde_1, \ldots, \bde_k)$ esimates $(\bmu_1, \ldots, \bmu_k)$.

A typical example of Model~(\ref{eqn:model}) is a $k$-sample problem.
For $i \in \{1, \ldots, k\}$, a random sample of size $n_i$ is drawn from the $i$th population, say 
$\X_{i1}, \ldots, \X_{i n_i} \sim \Nc_p(\bmu_i, \si^2 \V_{i,0})$,
where $\V_{i,0}$ is a known matrix.
In this case, $\X_i$, $\V_i$, $S$, and $n$ in (\ref{eqn:model}) correspond to 
$$
\Xb_i =  \sum_{j=1}^{n_i}\X_{ij}/n_i, \quad {\V_{i,0}}/{n_i}, \quad \tr \left\{\sum_{i=1}^k\V_{i,0}^{-1}\sum_{j=1}^{n_i}(\X_{ij}-\Xb_i)(\X_{ij}-\Xb_i)^\top \right\}, \quad \sum_{i=1}^k(n_i-1)p,
$$ 
respectively. Another example of (\ref{eqn:model}) is $k$ linear regression models such that
$\y_1 = \Z_1\bbe_1 + \bep_1,  \ldots, \y_k = \Z_k\bbe_k + \bep_k$, 
where $\y_i$ is an $n_i \times 1$ vector of observations, $\bbe_i$ is a $p \times 1$ vector of regression coefficients, $\Z_i$ is an $n_i\times p$ matrix of explanatory variables and $\bep_i$ is an $n_i \times 1$ vector having $\Nc_{n_i}(\zero, \si^2\I_{n_i})$. In this case, $\X_i$, $\V_i$, $S$ and $n$ in (\ref{eqn:model}) respectively correspond to
$$
\bbeh_i=(\Z_i^\top\Z_i)^{-1}\Z_i^\top\y_i, \quad (\Z_i^\top\Z_i)^{-1}, \quad 
\sum_{i=1}^k \Vert\y_i-\Z_i\bbeh_i\Vert_{\I_{n_i}}^2, \quad \sum_{i=1}^k (n_i-p). 
$$

In this paper, we consider the case where the means are equal or nearly so. For example, suppose that  Laboratories $L_1, \ldots, L_k$ use similar instruments to measure several common characteristics.
One  may then suspect that the $k$ population mean vectors are nearly equal.
In such a situation, rather than estimating each mean by the corresponding sample mean separately, it might be preferable to use a pooled mean estimator. A classical way to address this problem is to compute an estimator which uses the pooled mean estimator upon acceptance of the null hypothesis $\mathcal{H}_0 : \bmu_1=\cdots = \bmu_k$, and to use separate sample means upon rejection of $\mathcal{H}_0$.
However, the preliminary test estimators are non-smooth and do not necessarily improve on sample means.
As an alternative method, we consider Bayesian methods including hierarchical Bayes and empirical Bayes estimators.
  
A Bayesian approach to this problem uses the prior distribution of $\bmu_i$ assuming a multivariate normal distribution with common mean $\bnu$ and common variance $\tau^2$, where $\tau^2$ is assumed to have a hierarchical prior distribution.
The mean $\bnu$ corresponds to $\mathcal{H}_0 : \bmu_1=\cdots = \bmu_k=\bnu$.
Ghosh and Sinha \cite{gs1988} considered this prior distribution in the framework of estimating the single mean $\bmu_1$ in the two-sample problem and showed that the resulting empirical and hierarchical Bayes estimators have reasonable forms of shrinking $\X_1$ towards the pooled estimator of $\bnu$.
They also derived conditions for those shrinkage estimators to be minimax.
Shrinkage estimators with minimaxity have been studied since Stein \cite{stein1956} established the inadmissibility of the standard estimator $\X_1$.
For papers related to the context of this paper, see \cite{em1976, js1961, perron1993, stein1981, straw1971, straw1973, s1996}.

The Bayesian approach suggested here consists of two kinds of shrinkage.
One is to shrink the $k$ sample means toward a pooled mean estimator, and the other is to shrink the $p$-dimensional pooled estimator toward a constant like zero, which was proposed in~\cite{perron1993}.
The former shrinkage is treated in Section~\ref{sec:2.1}, and a class of minimax estimators is derived.
Out of the class, we develop hierarchical and empirical Bayes estimators with minimaxity.
The latter shrinkage is discussed in Section~\ref{sec:2.2}, and we show that all the estimators derived by the former shrinkage can be further improved by using the latter shrinkage.

In Section~\ref{sec:BM}, we find a hierarchical Bayes and minimax estimator, thus an admissible and minimax estimator using both kinds of shrinkage when $\Q_1=\V_1^{-1}, \ldots, \Q_k=\V_k^{-1}$.
The Bayes estimator, which consists of the two kinds of shrinkage, has a new form which we cannot handle with conventional techniques considered in Section~\ref{sec:2}.
Thus, we need to extend the classes of minimax estimators in order to treat such a new type of shrinkage estimator.
The minimaxity of the hierarchical Bayes estimator is successfully established.

In Section~\ref{sec:sim}, we investigate through simulations the performance of several hierarchical Bayes and empirical Bayes estimators and the preliminary test estimators.
The results show that the hierarchical Bayes and empirical Bayes estimators for the two kinds of prior distributions perform well. All proofs are relegated to the Appendix.

\section{Improvement by two kinds of shrinkage\label{sec:2}}

\subsection{Shrinkage toward the pooled estimator\label{sec:2.1}}

Let $\eta=1/\si^{2}$, $\bMu=(\bmu_1,\dots,\bmu_k)$ and $q=\min \{\Ch_{{\rm min}}(\V_1\Q_1), \ldots, \Ch_{{\rm min}}(\V_k\Q_k) \}$. For $\la \in (0, q)$, we begin by assuming the prior distribution 
\begin{equation}
\vec(\bMu)\mid\bnu,\eta,\la \sim \Nc_{pk}\left[\J\bnu,\;\;\blockdiag 
\left\{\frac{(\V_1\Q_1-\la \I_p)\V_1}{\la\eta},\dots,\frac{(\V_k\Q_k-\la \I_p)\V_k}{\la\eta}\right\} \right], \quad 
\bnu \sim \mathcal{U}(\Re^p),
\label{eqn:prior1}
\end{equation}
where $\J$ is the $pk\times p$ matrix $\J=(\I_p,\dots,\I_p)^\top$ and
$\mathcal{U}(\Re^p)$ denotes the improper uniform distribution over $\Re^p$.
The posterior distribution of $\bMu$ given $\X=(\X_1, \ldots, \X_k)$ and $\bnu$ is
$$
\vec(\bMu) \mid \X, \bnu, \eta, \la \sim \Nc_{pk} \left[ \vec\big\{\bMuh^*(\eta,\la,\bnu)\big\}, \blockdiag
\left( {\V_1-\la\Q_1^{-1} \over \eta},\dots,{\V_k-\la\Q_k^{-1} \over \eta} \right) \right]
$$
and the marginal distribution of $\X$ given $\bnu$ is
\begin{equation}
\vec(\X) \mid \bnu, \eta,\la \sim \Nc_{pk} \left[ \J\bnu, \blockdiag
\left( {\V_1\Q_1\V_1\over\la\eta},\dots,{\V_k\Q_k\V_k\over\la\eta} \right) \right],
\label{eqn:posterior1}
\end{equation}
where $\bMuh^*(\eta,\la,\bnu)=\big(\bmuh_1^*(\eta,\la,\bnu),\dots,\bmuh_k^*(\eta,\la,\bnu)\big)$ and 
\begin{equation}
\bmuh_i^*(\eta,\la,\bnu)= \X_i - \la\Q_i^{-1}\V_i^{-1}(\X_i-\bnu)=\X_i - \la\eta  \Q_i^{-1}\V_i^{-1}(\X_i-\bnu)/\eta.
\label{eqn:Bayes1}
\end{equation}
Also, the posterior distribution of $\bnu$ given $\X_1, \ldots, \X_k$ and the marginal density of $\X_1, \ldots, \X_k$ are
\begin{equation}
\begin{split}
\bnu \mid \X_1, \ldots, \X_k, \eta, \la \sim& \Nc_p[ \bnuh,\A / (\la\eta)],\\
f_\pi (\x_1, \ldots, \x_k \mid \eta, \la) \propto& (\la\eta)^{p(k-1)/2}\exp
\left (- \la\eta \sum_{i=1}^k\Vert\x_i-\bnuh\Vert_{\V_i^{-1}\Q_i^{-1}\V_i^{-1}}^2/2  \right),
\end{split}
\label{eqn:posterior2}
\end{equation}
where $\A = (\V_1^{-1}\Q_1^{-1}\V_1^{-1}+ \cdots + \V_k^{-1}\Q_k^{-1}\V_k^{-1})^{-1}$ and $\bnuh = \V_i^{-1}\Q_i^{-1}\V_1^{-1}\X_1 + \cdots + \V_k^{-1}\Q_i^{-1}\V_k^{-1}\X_i$.
Then, the Bayes estimator of $\bmu_i$ with respect to the quadratic loss (\ref{eqn:loss}) is derived from (\ref{eqn:Bayes1}) as
\begin{equation}
\bmuh_i^B(\eta, \la)=\bmuh_i^*(\eta, \la,\bnuh)=\X_i - \la\Q_i^{-1}\V_i^{-1}(\X_i-\bnuh).
\label{eqn:Bayes2}
\end{equation}

Because $\la$ and $\eta$ are unknown, we set
$$
F = \frac{1}{S} \sum_{i=1}^k \Vert \X_i-\bnuh\Vert_{\V_i^{-1}\Q_i^{-1}\V_i^{-1}}^2
$$
and estimate $\la\eta$ by $\{p(k-1)-2\}/(FS)$ from the marginal likelihood in (\ref{eqn:posterior2}) and $\si^2$ by $\sih^2=S/(n+2)$. Then, the resulting empirical Bayes estimator is
\begin{equation}
\bmuh_i^{EB1} = \X_i - \min \left[  {\{p(k-1)-2\}/\{(n+2)F\}}, 1 \right]\Q_i^{-1}\V_i^{-1}(\X_i-\bnuh).
\label{eqn:EB1}
\end{equation}
Motivated from the empirical Bayes estimator, we consider the class of estimators
\begin{equation}
\bmuh_i^S(\phi_0) = \X_i - \phi_0(F)  \Q_i^{-1}\V_i^{-1}(\X_i-\bnuh)/F,
\label{eqn:class1}
\end{equation}
where $\phi_0(F)$ is an absolutely continuous function.
We can provide conditions on $\phi_0(F)$ for minimaxity of $\bmuh_i^S(\phi_0)$. Theorem~\ref{thm:1} is a corollary of Proposition~\ref{prp:1} in Section~\ref{sec:BM}.

\begin{thm}
\label{thm:1}
The estimator $(\bmuh_1^S(\phi_0), \ldots, \bmuh_k^S(\phi_0))$ is minimax relative  to the quadratic loss $(\ref{eqn:loss})$ if
{\rm (a)} $\phi_0(F)$ is non-decreasing in $F$ and {\rm (b)} $0<\phi_0(F)\leq 2 \{p(k-1)-2\}/(n+2)$.
\end{thm}

It follows from Theorem~\ref{thm:1} that the empirical Bayes estimator $(\bmuh_1^{EB1}, \ldots, \bmuh_k^{EB1})$ is minimax. Another minimax estimator is the hierarchical Bayes estimator.
In addition to the prior distribution (\ref{eqn:prior1}), we assume that
$$
\pi(\la) \propto \;\la^{a -1} \mathbf{1}_{(0<\la<q)}, \quad
\pi(\eta) \propto \;\eta^{-c } \mathbf{1}_{(\eta>0)},
$$
where $a $ and $c $ are constants, $q=\min \{ \Ch_{{\rm min}}(\V_1 \Q_1), \ldots, \Ch_{{\rm min}}(\V_k \Q_k)\}$ and $\mathbf{1}$ denotes the indicator function.
It is noted that $\pi(\la)$ is proper for $a>0$, and $\pi(\eta)$ is improper for any $c$.
From (\ref{eqn:posterior2}), the posterior distribution of $(\la, \eta)$ given $\X_1, \ldots, \X_k, S$ is
$$
\pi(\la, \eta \mid \x_1, \ldots, \x_k,S)
\propto  \la^{p(k-1)/2+a -1}\eta^{\{n+p(k-1)\}/2-c}\exp \left( - {\la\eta }\sum_{i=1}^k\Vert\x_i-\bnuh\Vert_{\V_i^{-1}\Q_i^{-1}\V_i^{-1}}^2/2 -\eta S/2 \right),
$$
where $\la\in (0,q)$ and $\eta>0 $.
Then from (\ref{eqn:Bayes2}), the hierarchical Bayes estimator of $\bmu_i$ relative to the quadratic loss (\ref{eqn:loss}) is written as
\begin{align}
\bmuh_i^{HB1} =& \X_i -  {\rm E}  (\la \mid \X_1,  \ldots, \X_k, S)\Q_i^{-1}\V_i^{-1}(\X_i - \bnuh),
\label{eqn:HB1}
\end{align}
where
\begin{align*}
 {\rm E}  (\la \mid \X_1,  \ldots, \X_k, S)={\int_0^q \int_0^\infty \la^{p(k-1)/2+a }\eta^{\{n+p(k-1)\}/2-c }\exp\{ - {(\eta S/2)}(\la F+1)\}d\eta d\la
\over
\int_0^q \int_0^\infty \la^{p(k-1)/2+a -1}\eta^{\{n+p(k-1)\}/2-c }\exp\{ - {(\eta S/2)}(\la F+1)\}d\eta d\la}.
\end{align*}
Making the transformations $x=F\la$ and $v=S\eta$ gives the expression
${\rm E}  (\la \mid \X_1,  \ldots, \X_k, S)= {\phi^{HB1}(F) /F}$,
where
$$
\phi^{HB1}(F)=
{\int_0^{qF} x^{p(k-1)/2+a }/ (1+x)^{\{n+p(k-1)\}/2+1-c } dx
\over
\int_0^{qF} x^{p(k-1)/2+a -1}/(1+x)^{\{n+p(k-1)\}/2+1-c }dx}\\
$$
We can derive a condition on $a $ and $c $ for minimaxity of the hierarchical Bayes estimator.

\begin{thm}
\label{thm:HB1}
The hierarchical Bayes estimator $(\ref{eqn:HB1})$ is minimax if $p(k-1)\geq 3$, $a +c < (n+2)/2$ and
\begin{equation}
\{2p(k-1)+n-2\}a +2\{p(k-1)-2\}c  \leq (n-2)p(k-1)/2 - 2n.
\label{eqn:HB1c}
\end{equation}
\end{thm}

\noindent
\textit{Proof.}
It suffices to show that $\phi^{HB1}(F)$ satisfies conditions (a) and (b) in Theorem~\ref{thm:1}. The numerator of the derivative of $\phi^{HB1}$ with respect to $F$ is proportional to
\begin{align*}
   \frac{(qF)^{p(k-1)/2+a-1}}{(1+qF)^{\{n+p(k-1)\}/2+1-c }}\int_{0}^{qF}\frac{x^{p(k-1)/2+a-1}(qF-x)}{(1+x)^{\{n+p(k-1)\}/2+1-c}}\, dx,
\end{align*}
which establishes condition (a). To check condition (b), we note that
$$
\int_0^\infty {x^\ell  (1+x)^{-m}}dx = \int_0^1 w^\ell (1-w)^{m-\ell-2}dw
=B(\ell+1, m-\ell-1),
$$
for the beta function $B$.
Using condition (a), we see that
\begin{align*}
\phi^{HB1}(F)\leq \lim_{F\to\infty} \phi^{HB1}(F)
& =
{\int_0^\infty x^{p(k-1)/2+a }/ (1+x)^{\{n+p(k-1)\}/2+1-c } dx
\over
\int_0^\infty x^{p(k-1)/2+a -1}/(1+x)^{\{n+p(k-1)\}/2+1-c }dx} 
= {p(k-1)+2a  \over n- 2(a +c )}.
\end{align*}
Thus, condition (b) is satisfied if $a +c < n/2$ and if 
$$
\{p(k-1)+2a\}/\{n- 2(a +c )\} \leq 2 \{p(k-1)-2\}/(n+2),
$$
which is equivalently rewritten by (\ref{eqn:HB1c}).
Hence, the minimaxity of the hierarchical Bayes estimator is established.
\hfill$\Box$ 

\begin{remark}{\rm
When $\bnu=\zero$, the prior distribution (\ref{eqn:prior1}) was treated in~\cite{berger1976,zrnan2013}.
There exists a nonsingular matrix $\B_i$ such that $\B_i\V_i\B_i^\top=\I_p$ and $\B_i\Q_i^{-1}\B_i^\top=\bDe_i$ for a diagonal matrix $\bDe_i$.
This reduces the problem to the simple model with $\V_i=\I_p$ and $\Q_i=\bDe_i^{-1}$.
When we assume the uniform prior for $\bnu$ in (\ref{eqn:prior1}), however, we cannot reduce the problem to such a simple situation.
}
\end{remark}

\begin{remark}{\rm
When $\eta$ and $\la$ are known, assuming the uniform distribution for $\bnu$ in (\ref{eqn:prior1}) yields the Bayes estimator $\bmuh_i^*(\eta, \la,\bnuh)$.
It is interesting to note that this Bayes estimator is also interpreted as an empirical Bayes estimator, because when $\bnu$ is an unknown parameter, $\bnuh$ is the maximum likelihood estimator (MLE) of $\bnu$  in the marginal likelihood 
$$
\vec(\X)\sim \Nc_{pk}   [\J\bnu, \;\;(\la\eta)^{-1}\blockdiag ({\V_1\Q_1\V_1},\dots,{\V_k\Q_k\V_k} ) ].
$$
Thus, the estimator $\bmuh_i^*(\eta, \la,\bnuh)$ is a Bayes empirical Bayes estimator.
}
\end{remark}

\subsection{Shrinkage of the pooled estimator}
\label{sec:2.2}

The empirical Bayes estimator $\bmuh_i^{EB1}$ and the hierarchical Bayes estimator $\bmuh_i^{HB1}$ shrink $k$ individual estimators $\X_i$ toward the common mean $\bnuh$. 
Because $\bnuh$ is a $p$-dimensional estimator against the uniform prior over $\Re^p$, we can consider to shrink $\bnuh$ for further improvement. Such a further improvement was studied in~\cite{perron1993}. We here assume the prior distribution 
\begin{align}
\vec(\bMu)\mid\bnu,\eta,\la \sim& \Nc_{pk} \left[\J\bnu,\;\;\blockdiag
\{ {(\V_1\Q_1-\la I_p)\V_1}/{(\la\eta)},\dots, {(\V_k\Q_k-\la I_p)\V_k}/{(\la\eta)} \} \right],\non\\
\bnu \mid \eta,\la,\xi  \sim& \Nc_p [\zero, (\la-\xi)\A/(\la\xi\eta) ],
\label{eqn:prior3}
\end{align} 
where $\xi$ and $\la$ satisfy the constraint $0<\xi<\la<q$ for $q=\min\{\Ch_{{\rm min}}(\V_1 \Q_1), \ldots, \Ch_{{\rm min}}(\V_k \Q_k)\}$. 
It is noted that the normal prior is considered for $\bnu$ instead of the uniform prior in (\ref{eqn:prior1}). 
Then the posterior distribution of $\bMu$ given $\X=(\X_1, \ldots, \X_k)$ and $\bnu$, and the marginal distribution of $\X$ given $\bnu$ are given by (\ref{eqn:posterior1}) and (\ref{eqn:Bayes1}), and the posterior distribution of $\bnu$ given $\X_1, \ldots, \X_k$ and the marginal density of $\X_1, \ldots, \X_k$ are
\begin{equation}
\begin{split}
\bnu \mid \X_1, \ldots, \X_k, \eta, \la,\xi \sim& \Nc_p [ \left(1- {\xi}/{\la}\, \right)\bnuh, {(\la-\xi)}\A/{(\la^2\eta)} ],\\
f_\pi(\x_1,\ldots, \x_k \mid \la, \xi, \si^2) \propto&
 ({\la\eta})^{p(k-1)/2}\exp \left( -{\la\eta } \sum_{i=1}^k\Vert\x_i-\bnuh\Vert_{\V_i^{-1}\Q_i^{-1}\V_i^{-1}}^2/2  \right)
\times  \left({\xi\eta}\right)^{p/2} \exp ( -{\xi\eta } \Vert\bnuh\Vert_{\A^{-1}}^2/2  ).
\end{split}
\label{eqn:marginal3}
\end{equation}
Thus the Bayes estimator of $\bmu_i$ with respect to the prior in (\ref{eqn:prior3}) is given by 
\begin{equation}
\bmuh_i^B(\eta, \la,\xi)=
\X_i - \la\Q_i^{-1}\V_i^{-1}(\X_i-\bnuh)-\xi\Q_i^{-1}\V_i^{-1}\bnuh.
\label{eqn:Bayes3}
\end{equation}

Since $\la$, $\xi$ and $\eta$ are unknown, we now set
$$
F = \frac{1}{S} \sum_{i=1}^k\Vert\X_i-\bnuh \Vert^2_{\V_i^{-1}\Q_i^{-1}\V_i^{-1}}, 
$$
and estimate $\la\eta$ and $\xi\eta$ by $\{p(k-1)-2\}(FS)$ and $(p-2)/\Vert\bnuh\Vert_{\A^{-1}}^2$, respectively. When $\si^2$ $(=1/\eta)$ is estimated by $\sih^2=S/(n+2)$, the resulting empirical Bayes estimator is
\begin{equation}
\bmuh_i^{EB2} = \X_i - \min\left[ {\{p(k-1)-2\}/\{(n+2)F\}}, 1\right] \Q_i^{-1}\V_i^{-1}(\X_i-\bnuh)
 -\min\left[ {(p-2)/\{(n+2)G\}}, 1\right] \Q_i^{-1}\V_i^{-1}\bnuh,
\label{eqn:EB2}
\end{equation}
where $G=\Vert\bnuh\Vert_{\A^{-1}}^2/S$.
The estimator $\bmuh_i^{EB2}$ enjoys shrinking $\X_i$ toward $\bnuh$ and shrinking $\bnuh$ toward zero.

Motivated from the empirical Bayes estimator, we consider the class of double shrinkage estimators
\begin{equation}
\bmuh_i^S(\phi_0, \psi_0) = \X_i - {\phi_0(F)} \Q_i^{-1}\V_i^{-1}(\X_i-\bnuh)/F - {\psi_0(G)} \Q_i^{-1}\V_i^{-1}\bnuh/G,
\label{eqn:class2}
\end{equation}
where $\phi_0(F)$ and $\psi_0(G)$ are absolutely continuous functions.
We can provide conditions on $\psi_0(G)$ for improving on $\bmuh_i^S(\phi_0)$ in (\ref{eqn:class1}), which is stated in the following theorem. Theorem~\ref{thm:2} is a corollary of Proposition~\ref{prp:2} in the subsequent section.

\begin{thm}
\label{thm:2}
The estimator $(\bmuh_1^S(\phi_0), \ldots, \bmuh_k^S(\phi_0))$ in $(\ref{eqn:class1})$ is improved upon by the double shrinkage estimator $(\bmuh_1^S(\phi_0, \psi_0)$, $\ldots$, $\bmuh_k^S(\phi_0,\psi_0))$ in $(\ref{eqn:class2})$ relative  to the quadratic loss $(\ref{eqn:loss})$ if 
{\rm (a)} $\psi_0(G)$ is non-decreasing in $G$, and
{\rm (b)} $0<\psi_0(G)\leq 2 (p-2)/(n+2)$.
\end{thm}

From Theorem~\ref{thm:2}, the empirical Bayes estimator $(\bmuh_1^{EB1}, \ldots, \bmuh_k^{EB1})$ in (\ref{eqn:EB1}) is improved upon by $(\bmuh_1^{EB2}, \ldots, \bmuh_k^{EB2})$ in (\ref{eqn:EB2}). 
Using Theorem~\ref{thm:2}, we provide a shrinkage estimator improving on the preliminary-test (PT) estimator.
We here consider the case where $\Q_i=\V_i^{-1}$ for all $i\in\{1, \ldots, k\}$. 
With the new definition of $F$, the statistic $nF/\{p(k-1) \}$ has the Fisher--Snedecor distribution with $(p(k-1),n)$ degrees of freedom under the hypothesis $\mathcal{H}_0$. 
Then the PT estimator based on the statistic $F$ for testing the hypothesis is
\begin{equation}
\bmuh_i^{PT}= \X_i - \mathbf{1}_{[F\leq \{p(k-1)/n\}F_{p(k-1), n, \al}]}(\X_i -\bnuh)  =
\left\{\begin{array}{ll} \X_i & \ \text{if}\ F>\{p(k-1)/n\}F_{p(k-1), n, \al}, \\
\bnuh & \ \text{otherwise,}
\end{array}\right.
\label{eqn:PT}
\end{equation}
where $F_{p(k-1), n, \al}$ is the upper $\al$ point of the Fisher--Snedecor distribution with $(p(k-1), n)$ degrees of freedom.
Because the PT estimator belongs to the class (\ref{eqn:class1}), Theorem~\ref{thm:2} implies that the estimator $(\bmuh_1^{PT}, \ldots, \bmuh_k^{PT})$ can be improved upon by $(\bmuh_1^{PT*}, \ldots, \bmuh_k^{PT*})$ with
\begin{equation}
\bmuh_i^{PT*} = \bmuh_i^{PT} - \min \left[ (p-2)/\{( n+2) G\}, 1 \right] \bnuh.
\label{eqn:PTI}
\end{equation}

\section{Bayes and minimax estimation\label{sec:BM}}

\subsection{Extension of the class of minimax shrinkage estimators}

In this section, we derive the hierarchical Bayes and minimax estimators.
As shown in (\ref{eqn:HB2}), the hierarchical Bayes estimator with double shrinkage does not belong to the class (\ref{eqn:class2}).
To show minimaxity of the hierarchical Bayes estimators, we begin by extending the class of estimators (\ref{eqn:class2}) to  the following class with more general form
\begin{equation}
\bmuh_i(\phi, \psi) = 
\X_i - \phi(F,G,S)  \Q_i^{-1}\V_i^{-1}(\X_i - \bnuh)/F - \psi(F,G, S) \Q_i^{-1}\V_i^{-1}\bnuh/G,
\label{eqn:GE}
\end{equation}
where $\phi(F,G,S)$ and $\psi(F,G,S)$ are absolutely continuous functions.
The class (\ref{eqn:class2}) is a special case of (\ref{eqn:GE}), and Theorems \ref{thm:1} and \ref{thm:2} can be extended in the propositions given below.

\begin{thm}
\label{thm:UER}
An unbiased estimator of the risk function of $(\ref{eqn:GE})$ relative to the loss $(\ref{eqn:loss})$ is 
\begin{align}
\text{UER}(\phi,\psi)=&
\sum_{i=1}^k\tr(\V_i\Q_i)  - {[2\{p(k-1)-2\}-(n+2)\phi]\phi/ F} - 4 \phi_F - {4\phi }(F\phi_F+G\phi_G)/F+{4S\phi\phi_S/F}\non\\
& - {\{2(p-2)-(n+2)\psi\}\psi/ G} - 4 \phi_G - {4\psi}(F\psi_F+G\psi_G)/G +{4S\psi\psi_S/ G},
\label{eqn:UER}
\end{align}
where $\phi_F=\partial \phi(F,G,S)/\partial F$, $\phi_G=\partial \phi(F,G,S)/\partial G$, $\phi_S=\partial \phi(F,G,S)/\partial S$, and $\psi_F$, $\psi_G$ and $\psi_S$ are defined similarly.
\end{thm}

The proof is given in the Appendix.
Considering the case with $\psi(F,G,S)=0$ in (\ref{eqn:UER}), we can get a generalization of Theorem~\ref{thm:1}.

\begin{prp}
\label{prp:1}
The estimator $(\bmuh_1(\phi), \ldots, \bmuh_k(\phi))$ with the form $\bmuh_i(\phi) = \X_i - \{\phi(F,G,S)/F\}(\X_i - \bnuh)$ is minimax relative  to the quadratic loss $(\ref{eqn:loss})$ if 
{\rm (a)} $\phi(F,G,S)$ is non-decreasing in $F$ and  $G$,  and non-increasing in $S$; and 
{\rm (b)} $0<\phi(F,G,S)\leq 2 \{p(k-1)-2\}/(n+2)$.
\end{prp}

The unbiased risk estimator (\ref{eqn:UER}) also enables us to get a generalization of Theorem~\ref{thm:2}.

\begin{prp}
\label{prp:2}
The estimator $(\bmuh_1(\phi), \ldots, \bmuh_k(\phi))$ given in Proposition~\ref{prp:1} is improved upon by the double shrinkage estimator $(\bmuh_1(\phi, \psi), \ldots, \bmuh_k(\phi,\psi))$ in $(\ref{eqn:GE})$ relative  to the quadratic loss $(\ref{eqn:loss})$ if 
{\rm (a)} $\psi(F,G,S)$ is non-decreasing in $F$ and $G$, and  non-increasing in $S$; and 
{\rm (b)} $0<\psi(F,G,S)\leq 2 (p-2)/(n+2)$.
\end{prp}

\subsection{Hierarchical Bayes minimax estimators}

In this section, we concentrate attention on the case of $\Q_1=\V_1^{-1}, \ldots, \Q_k=\V_k^{-1}$.
Although a natural prior is a hierarchical prior based on  (\ref{eqn:prior3}), the minimaxity of the resulting Bayes estimator is difficult to show, because the integrals with respect to $\xi$ and $\la$ are taken over the constraint $\xi<\la$.
Instead of (\ref{eqn:prior3}), we here consider the following prior distribution for $\bMu=(\bmu_1,\dots,\bmu_k)$:
\begin{equation}
    \vec(\bMu)\mid \eta,\la,\xi\sim \Nc_{pk} \left[0,\;\;\eta^{-1} \left\{ {\la\over 1-\la} \, \big (\bSi^{-1}-\bSi^{-1}\J\A\J^\top\bSi^{-1}\big)+{\xi\over1-\xi  }\, \bSi^{-1}\J\A\J^\top\bSi^{-1} \right\}^{-1} \right],
	\label{eqn:modified prior}
\end{equation}
 where $0<\la<1$, $0<\xi<1$, $\bSi=\blockdiag(\V_1,\dots,\V_k)$ and $\A=\big(\V_1^{-1} + \cdots + \V_k^{-1}\big)^{-1}$. 
The covariance matrix in the prior (\ref{eqn:modified prior}) is positive definite for $0<\la, \xi<1$. 
It is noted that when $\Q_i=\V_i^{-1}$, integrating out the prior (\ref{eqn:prior3}) with respect to $\bnu$ leads to the same density function as in (\ref{eqn:modified prior}) with the constraint $\xi<\la$.
The prior distribution (\ref{eqn:modified prior}) does not have to assume this constraint.

For $\la$, $\xi$ and $\eta$, we assume the following second-stage priors:
\begin{equation}
\pi(\la) \propto \;\la^{a -1}\mathbf{1}_{(0<\la< 1)},\quad
\pi(\xi) \propto \;\xi^{b -1}\mathbf{1}_{(0<\xi< 1)},\quad
\pi(\eta) \propto \;\eta^{-c}\mathbf{1}_{(\eta\geq L)}.
\label{eqn:prior4}
\end{equation}
Here $a $, $b$ and $c $ are constants and $L $ is a positive constant.
When $a>0$, $b>0$ and $c>1$, the prior distribution (\ref{eqn:prior4}) is proper. 
From (\ref{eqn:marginal3}), the posterior distribution of $(\la, \xi, \eta)$ given $\X_1, \ldots, \X_k, S$ is
$$
\pi(\xi, \la, \eta\mid \x_1, \ldots, \x_k, S)
\propto
\la^{p(k-1)/2+a-1} \xi^{p/2+b-1} \eta^{(n+pk)/2-c}\exp\left\{ -{S \eta} (\la F+ \xi G + 1)/2\right\},
$$
where $0<\xi<q$, $0<\la<q$ and $\eta\geq L $. From (\ref{eqn:Bayes3}), the hierarchical Bayes estimator of $\bmu_i$ is written as
\begin{align}
\bmuh_i^{HB2} = 
\X_i -  {\rm E}  (\la \mid \X, S)(\X_i - \bnuh)
-  {\rm E}  (\xi \mid \X, S)\bnuh,
\label{eqn:HB20}
\end{align}
where $\X=(\X_1, \ldots, \X_k)$,
\begin{align*}
 {\rm E}  (\la \mid \X, S)={\int_0^1\int_0^1 \int_L^\infty \la^{p(k-1)/2+a }\xi^{p/2+b-1}\eta^{(n+pk)/2-c }\exp\{ - {\eta Sr 2}(\la F+\xi G +1)/2\}d\eta d\xi d\la
\over
\int_0^1 \int_0^1 \int_L ^\infty \la^{p(k-1)/2+a-1 }\xi^{p/2+b-1}\eta^{(n+pk)/2-c }\exp\{ - {\eta S }(\la F+\xi G + 1)/2\}d\eta d\xi d\la}.
\end{align*}
Making the transformations $x=F\la$, $y=G\xi$ and $v=S\eta$ gives the expression
$ {\rm E}  (\la \mid \X, S)= {\phi^{HB2}(F,G,S)/ F}$,
where
$$
\phi^{HB2}(F,G,S)=
{\int_0^{F} \int_0^{G} \int_{L  S}^\infty x^{p(k-1)/2+a }y^{p/2+b-1}v^{(n+pk)/2-c }\exp\{ - v(x+y+1)/2\}dv dy dx
\over
\int_0^{F} \int_0^{G} \int_{L  S}^\infty x^{p(k-1)/2+a-1 }y^{p/2+b-1}v^{(n+pk)/2-c }\exp\{ - v(x+y+1)/2\}dv dy dx}.
$$
Similarly,
$ {\rm E}  (\xi \mid \X, S)= {\psi^{HB2}(F,G, S) /G}$,
where
$$
\psi^{HB2}(F,G,S)=
{\int_0^{F} \int_0^{G} \int_{L  S}^\infty x^{p(k-1)/2+a-1 }y^{p/2+b}v^{(n+pk)/2-c }\exp\{ - v(x+y+1)/2\}dv dy dx
\over
\int_0^{F} \int_0^{G} \int_{L  S}^\infty x^{p(k-1)/2+a-1 }y^{p/2+b-1}v^{(n+pk)/2-c }\exp\{ - v(x+y+1)/2\}dv dy dx}.
$$
Thus, the hierarchical Bayes estimator in (\ref{eqn:HB20}) is expressed as
\begin{equation}
\bmuh_i^{HB2} = 
\X_i - {\phi^{HB2}(F,G,S) }(\X_i - \bnuh)/F - {\psi^{HB2}(F,G, S) }\bnuh/G.
\label{eqn:HB2}
\end{equation}

This estimator does not belong to the class (\ref{eqn:class2}), but belongs to the class (\ref{eqn:GE}).
For minimaxity of the hierarchical Bayes estimators, we check the conditions given in Propositions~\ref{prp:1}--\ref{prp:2}.
Combining the condition for the minimaxity and the condition for the proper prior, we can get a condition for the Bayes estimator to be admissible and minimax.

\begin{thm}
\label{thm:HB2}
The hierarchical Bayes estimator $(\ref{eqn:HB2})$ is minimax if $p(k-1)\geq 3$, $p\geq 3$ and if  $a +b+c < n/2$ and
\begin{equation}
\{2p(k-1)+n-2\}a +2\{p(k-1)-2\}(b+c)  \leq p(k-1)(n-2)/2-2n,
\label{eqn:HB2c}
\end{equation}
\begin{equation}
(2p+n-2)b +2(p-2)(a+c)  \leq p(n-2)/2-2n.
\label{eqn:HB2cc}
\end{equation}
Further, $(\ref{eqn:HB2})$ is a Bayes minimax estimator, namely admissible and minimax, provided that the constants $a$, $b$ and $c$ satisfy the above conditions for $a>0$, $b>0$ and $c>1$.
\end{thm}

The latter part of Theorem~\ref{thm:HB2} follows from the fact that the prior distribution (\ref{eqn:prior4}) is proper when $a >0$, $b >0$, $c>1$ and $L >0$.
Then, the Bayes estimator (\ref{eqn:HB2}) is admissible and minimax.
Since $k\geq 2$, we can find such constants if  $p\geq 5$ and $n>2(3p-4)/(p-4)$.  
This shows that $p$ is at least $5$.

\begin{remark}{\rm
Although the condition $L>0$ is imposed technically in order to guarantee the admissibility, we have no information on $L$. Furthermore, the Bayes estimator based on triple integrals is computationally hard to derive.
From a practical point of view, we set $L=0$.
In this case, the hierarchical prior (\ref{eqn:prior4}) is improper and $\bmuh_i^{HB2}$ is no longer guaranteed to be admissible. 
However, the minimaxity still holds, and the resulting generalized Bayes estimator is
\begin{equation}
\bmuh_i^{HB2}=\X_i - {\phi^{HB2}(F,G)}(\X_i - \bnuh)/F - {\psi^{HB2}(F,G)} \bnuh/G,
\label{eqn:HB2a}
\end{equation}
where $\phi^{HB2}(F,G,S)$ and $\psi^{HB2}(F,G,S)$ are expressed based on double integrals as
\begin{align*}
	&\phi^{HB2}(F,G)={\int_0^{F} \int_0^{G} x^{p(k-1)/2+a }y^{p/2+b-1} (x+y+1)^{-\{(n+pk)/2-c+1 \}} dy dx
		\over
		\int_0^{F} \int_0^{G}  x^{p(k-1)/2+a-1 }y^{p/2+b-1}(x+y+1)^{-\{(n+pk)/2-c+1 \}} dy dx},\\
	&\psi^{HB2}(F,G)={\int_0^{F} \int_0^{G} x^{p(k-1)/2+a-1 }y^{p/2+b}(x+y+1)^{-\{(n+pk)/2-c+1 \}} dy dx
		\over
		\int_0^{F} \int_0^{G}  x^{p(k-1)/2+a-1 }y^{p/2+b-1}(x+y+1)^{-\{(n+pk)/2-c+1 \}} dy dx}.
\end{align*}
In our simulation experiments given in the next section, we use this generalized Bayes estimator.
}\end{remark}

\section{Simulation study}
\label{sec:sim}

We investigate the numerical performances of the risk functions of the preliminary-test estimator and several empirical and hierarchical Bayes estimators through simulation. To compare the preliminary-test and Bayes estimators fairly, we employ the quadratic loss function $L(\bde,\bom)$ in (\ref{eqn:loss}) for $ \Q_1= \V_1^{-1}, \ldots, \Q_k= \V_k^{-1}$.
The estimators which we compare are the following eight:
\begin{itemize}
\item []
JS1: the sample-wise James--Stein estimator
$$
\bmuh_i^{JS1}=\X_i - {p-2 \over n+2}{S\over \Vert\X_i\Vert^2_{\V_i^{-1}}}\X_i,
$$
\item []
JS2: the overall James--Stein estimator for $\bMu=(\bmu_1,\dots,\bmu_k)$
$$
\widehat{\boldmath{\text{$\mathcal{M}$}}}^{JS2}=\X - {pk-2 \over n+2}{S\over \sum_{i=1}^k\Vert\X_i\Vert^2_{\V_i^{-1}}}\X,
$$
\item []
PT: the preliminary-test estimator $\bmuh_i^{PT}$ given in (\ref{eqn:PT}), 
\item []
PT*: the preliminary-test and shrinkage estimator $\bmuh_i^{PT*}$ in (\ref{eqn:PTI}), 
\item []
EB: the empirical Bayes estimator $\bmuh_i^{EB1}$ in (\ref{eqn:EB1}),
\item []
EB*: the improved empirical Bayes estimator $\bmuh_{EB2}$ in (\ref{eqn:EB2})
\item []
HB1: the hierarchical Bayes estimator $\bmuh_i^{HB1}$ in (\ref{eqn:HB1}),
\item []
HB2: the hierarchical Bayes estimator $\bmuh_i^{HB2}$ in (\ref{eqn:HB2a}).
\end{itemize}
  
\begin{table}[b!]
\caption{Values of PRIAL of estimators PT, PT*, EB1, EB2, HB1 and HB2}
\begin{center}
$
{\renewcommand\arraystretch{1.1}\small
\begin{array}{l@{\hspace{5mm}}
              r@{\hspace{2mm}}
              r@{\hspace{2mm}}
              r@{\hspace{2mm}}
              r@{\hspace{2mm}}
              r@{\hspace{2mm}}
              r@{\hspace{2mm}}
              r@{\hspace{2mm}}
              r
             }
\text{$(\bmu_1, \ldots, \bmu_5)$} &\text{JS1}&\text{JS2}&\text{PT}&\text{PT*}&\text{EB}&\text{EB*}&\text{HB1}&\text{HB2}\\
\hline
\text{$(\zero, \zero, \zero, \zero, \zero)$} 
&54.1 
&82.8
&74.1 
&87.2 
&70.3 
&83.4 
&69.7 
&82.6 
\\
\text{$(2, 2, 2, 2, 2)\otimes \j_5$} 
&9.0
&14.2
&74.1 
&74.4 
&70.3 
&70.6 
&69.7
&70.0 
\\
\hline

\text{$(-0.4, -0.2, 0.0, 0.2, 0.4)\otimes \j_5$} 
&48.3
&74.4
&64.1 
&75.3 
&63.6 
&74.9 
&64.7 
&76.5 
\\
\text{$(-1, -0.5, 0, 0.5, 1)\otimes \j_5$} 
&37.4
&48.2
&17.4 
&23.0 
&40.6 
&46.2 
&44.6 
&52.2 
\\
\text{$(-2, -1, 0, 1, 2)\otimes \j_5$} 
&25.0
&21.1
&-11.6
&-10.0
&17.5 
&19.1 
&17.5
&18.9
\\
\hline
\text{$(1.2, 1.4, 1.6, 1.8, 2.0)\otimes \j_5$} 
&12.7
&23.3
&64.1 
&64.6 
&63.6 
&64.2 
&64.7 
&65.2 
\\
\text{$(1, 1.5, 2, 2.5, 3)\otimes \j_5$} 
&9.5
&18.1
&17.4 
&17.9 
&40.6 
&41.1 
&44.6 
&44.9 
\\
\text{$(0, 1, 2, 3, 4)\otimes \j_5$} 
&18.9
&17.1
&-11.6 
&-10.8 
&17.5 
&18.3 
&17.5 
&17.6 
\\
\end{array}
}
$
\end{center}
\label{table:risk1}
\end{table}

The significance size  in the preliminary-test and related estimators $\bmuh_i^{PT}$ and $\bmuh_i^{PT*}$ is $\al=0.05$.
For the hierarchical Bayes estimators $\bmuh_i^{HB1}$ and $\bmuh_i^{HB2}$, we used the constants $a=b=c=0.1$ and $L=0$. 
In this simulation, we generate random numbers of $\X_1, \ldots, \X_k$ and $S$ based on Model~(\ref{eqn:model}) for $p=k=5$, $n=20$, $\si^2=2$ and $\V_i=(0.1\times i)\, \I_p$ for all $i\in\{1, \ldots, k\}$.
For the mean vectors $\bmu_i$, we treat the eight cases: 
\begin{multline*}(\bmu_1, \ldots, \bmu_5)  = (\zero, \zero, \zero, \zero, \zero), (2\j_5, 2\j_5, 2\j_5, 2\j_5, 2\j_5),\\
 (-0.4\j_5, -0.2\j_5, \zero, 0.2\j_5, 0.4\j_5), (-\j_5, -0.5\j_5, \zero, 0.5\j_5, \j_5), (-2\j_5, -\j_5, \zero, \j_5, 2\j_5),
\\
 (1.2\j_5, 1.4\j_5, 1.6\j_5, 1.8\j_5, 2.0\j_5), (\j_5, 1.5\j_5, 2\j_5, 2.5\j_5, 3\j_5),(\zero, \j_5, 2\j_5, 3\j_5, 4\j_5),  
\end{multline*}
where $\j_p=(1, \ldots, 1)^\top\in \Re^p$.
The first two are cases of equal means, the next three are cases where $\bmu_1 + \cdots + \bmu_5 =\zero$, and the last three are unbalanced cases.
  
For each estimator $\bde=(\bmuh_1, \ldots, \bmuh_5)$, based on 5000 replicates, we obtain an approximated value of the risk function $R(\bom, \bde)={\rm E} \{ L(\bde, \bom) \}$ for the loss function given in (\ref{eqn:loss}).
Table~\ref{table:risk1} reports the percentage relative improvement in average loss (PRIAL) of each estimator $\bde$ over $\bde^U=(\X_1, \ldots, \X_5)$, defined by
$$
{\rm PRIAL} = 100\{ R(\bom, \bde^U) - R(\bom, \bde)\}/R(\bom, \bde^U).
$$

Note from Theorem~\ref{thm:2} that PT and EB1 can be improved upon by PT* and EB2.
These results are confirmed by the simulation results in Table~\ref{table:risk1}.
Concerning the hierarchical Bayes estimators HB1 and HB2, Theorem~\ref{thm:HB2} shows that the estimator HB2 is minimax, while we could not demonstrate that HB2 dominates HB1.
However, Table~\ref{table:risk1} illustrates that HB2 is better than HB1.
  
Concerning PT, EB1 and HB1, their values of PRIAL are high under the hypothesis $\mathcal{H}_0 : \bmu_1=\cdots=\bmu_5$ but decrease when the means are far away from $\mathcal{H}_0$. 
The empirical Bayes estimator EB1 and the hierarchical Bayes estimator HB1 are comparable, and EB2 and HB2 are also comparable.
The simulation results in these setups illustrate that the hierarchical Bayes estimator HB2 and the empirical Bayes estimator EB2 have good performances.

\section{Proofs}

We here provide the proofs of Theorems \ref{thm:UER} and \ref{thm:HB2}.
Theorems \ref{thm:1} and \ref{thm:2} can be derived from Theorem~\ref{thm:UER} through Propositions \ref{prp:1} and \ref{prp:2}.

\subsection{Proof of Theorem~\ref{thm:UER}}

For the proof, Stein's identity \cite{stein1981} and the chi-square identity due to Efron and Morris \cite{em1976} are useful.
See also Bilodeau and Kariya \cite{bk1989} for a multivariate version of Stein's identity.

\begin{lem}
\label{lem:identity}
{\rm (i)} Assume that $\Y=(Y_1, \ldots, Y_p)^\top$ is a $p$-variate random vector having $\Nc_p(\bmu, \bSi)$ and that $\h$ is an absolutely continuous function from $\Re^p$ to $\Re^p$.
Then, we have Stein's identity, viz.
$$ 
{\rm E}  \{(\Y-\bmu)^\top \h(\Y)\}= {\rm E}  \big[ \tr\big\{ \bSi \bnabla_\Y \h(\Y)^\top\big\}\big],
$$
provided the expectations in both sides exist, where $\bnabla_\Y=(\partial /\partial Y_1, \ldots, \partial /\partial Y_p)^\top$.

{\rm (ii)} Assume that $S$ is a random variable such that $S/\si^2\sim \chi_{(n)}^2$ and that $g$ is an absolutely continuous function from $\Re$ to $\Re$.
Then, we have the chi-square identity, viz.
$$
 {\rm E}  \{S g(S)\} = \si^2  {\rm E}  \{n g(S) + 2Sg'(S)\},
$$
provided the expectations exist on both sides.
\end{lem}

The risk function is decomposed, say, as
$$
R\{ \bom, \bmuh(\phi, \psi)\}  = \frac{1}{\si^2}
\sum_{i=1}^k  {\rm E} \Big \{\Vert\bmuh_i(\phi,\psi)-\bmu_i\Vert^2_{\Q_i^{-1}} \Big\}  =  I_1 - 2 I_2 - 2I_3 + I_4, 
$$
where
$$
I_1 = \frac{1}{\si^2} \sum_{i=1}^k  {\rm E}  \Big\{\Vert\X_i-\bmu_i\Vert^2_{\Q_i^{-1}} \Big\} , \quad 
I_2 = {\phi\over \si^2 F} \sum_{i=1}^k  {\rm E} \{(\X_i-\bmu_i)^\top\V_i^{-1}(\X_i-\bnuh) \},
$$
$$
I_3 = {\psi\over \si^2 G} \sum_{i=1}^k  {\rm E} \{(\X_i-\bmu_i)^\top\V_i^{-1}\bnuh \}, \quad
I_4 =  {1\over \si^2}\sum_{i=1}^k  {\rm E} \left\{  \Vert{\phi }(\X_i-\bnuh)/F+{\psi}\bnuh/G \Vert^2_{\V_i^{-1}\Q_i^{-1}\V_i^{-1}} \right\}.
$$
with $\phi=\phi(F,G,S)$ and $\psi=\psi(F,G,S)$. It is easy to see $I_1 = \tr(\V_1\Q_1) + \cdots + \tr(\V_k\Q_k)$. Note that 
\begin{align*}
S\times F&=\sum_{i=1}^{k}\X_i^\top\V_i^{-1}\Q_i^{-1}\V_i^{-1}\X_i-
\left(\sum_{i=1}^{k}\V_i^{-1}\Q_i^{-1}\V_i^{-1}\X_i \right)^\top \A \left(\sum_{i=1}^{k}\V_i^{-1}\Q_i^{-1}\V_i^{-1}\X_i \right),\\
S\times G&= \left(\sum_{i=1}^{k}\V_i^{-1}\Q_i^{-1}\V_i^{-1}\X_i \right)^\top \A 
\left(\sum_{i=1}^{k}\V_i^{-1}\Q_i^{-1}\V_i^{-1}\X_i \right).
\end{align*}

Letting $\bnabla_i=\partial/\partial \X_i=(\partial /\partial X_{i1}, \ldots, \partial /\partial X_{ip})^\top$ for $\X_i=(X_{i1}, \ldots, X_{ip})^\top$,
we have $\bnabla_i F=2\V_i^{-1}\Q_i^{-1}\V_i^{-1}(\X_i-\bnuh)/S$ and $\bnabla_i G=2\V_i^{-1}\Q_i^{-1}\V_i^{-1}\bnuh/S$.
For $I_2$, it is observed that, by Stein's identity,
\begin{align*}
I_2& = \sum_{i=1}^k  {\rm E}  \Big[\tr [\bnabla_i  \{(\X_i-\bnuh)^\top {\phi(F,G,S)/F} \} ]\Big]
\\
& =
\sum_{i=1}^k  {\rm E}  \Big[ \tr\Big[  \{\bnabla_i (\X_i-\bnuh)^\top \} {\phi(F,G,S)/ F} +  \{\bnabla_i {\phi(F,G,S)/  F} \}(\X_i-\bnuh)^\top\Big]
\\
& =
\sum_{i=1}^k  {\rm E}  \Big[ \tr\Big[(\I_p-\V_i^{-1}\Q_i^{-1}\V_i^{-1}\A){\phi/F}\\
& \quad \quad +  (\X_i-\bnuh)^\top \{
{2\V_i^{-1}\Q_i^{-1}\V_i^{-1}(\X_i-\bnuh) } ( -{\phi/ F^2}+{\phi_F/ F} )/S+ {2 \V_i^{-1}\Q_i^{-1}\V_i^{-1}\bnuh } {\phi_G/(F}S) \}\Big]\Big]\\
& =
 {\rm E}   \left\{ p(k-1){\phi/ F} + 2 F ( -{\phi/ F^2}+{\phi_F/ F} )
+ \sum_{i=1}^k {(\X_i-\bnuh)^\top\V_i^{-1}\Q_i^{-1}\V_i^{-1}\bnuh } {\phi_G /(FS)}\right \}.
\end{align*}
Noting that $\sum_{i=1}^k (\X_i-\bnuh)^\top\V_i^{-1}\Q_i^{-1}\V_i^{-1}=\zero$, we get
\begin{equation}
I_2 =  {\rm E}  [ { \{p(k-1) -2\} \phi / F} + 2 \phi_F].
\label{eqn:I2}
\end{equation}
Similarly, we have
\begin{align*}
I_3 & = \sum_{i=1}^k  {\rm E}  \Big[\tr [\bnabla_i  \{\bnuh^\top {\psi(F,G,S)/ G} \} ]\Big]
\\
& =
\sum_{i=1}^k  {\rm E}  \Big[ \tr [(\V_i^{-1}\Q_i^{-1}\V_i^{-1}\A){\psi/ G} 
+  \bnuh^\top  \{
{2\V_i^{-1}\Q_i^{-1}\V_i^{-1}\bnuh }  ( -{\psi/ G^2}+{\psi_G/G} )/S + {2 \V_i^{-1}\Q_i^{-1}\V_i^{-1}(\X_i-\bnuh) }{\phi_F /(GS)}\} ]\Big].
\end{align*}
Since $\sum_{i=1}^k \bnuh^\top \V_i^{-1}\Q_i^{-1}\V_i^{-1}(\X_i-\bnuh)=\zero$, we get
\begin{equation}
I_3 =  {\rm E}  \{ (p-2)\psi / G + 2 \psi_G\}.
\label{eqn:I3}
\end{equation}

Concerning $I_4$, note that
\begin{align*}
I_4=&
{1\over \si^2} \, {\rm E}  \Big\{{S }\phi^2(F,G,S)/F\Big\} 
+ {1\over \si^2} \, {\rm E}  \Big\{{S}\psi^2(F,G,S)/G\Big\} 
+ {2\over \si^2} \, {\rm E}  \left\{\sum_{i=1}^k(\X_i-\bnuh)^\top\V_i^{-1}\Q_i^{-1}\V_i^{-1}\bnuh {\phi\psi/(FG)}\right\}
\\
=& {1\over \si^2} \, {\rm E} \{{S}\phi^2(F,G,S)/F \} 
+ {1\over \si^2} \, {\rm E} \{{S }\psi^2(F,G,S)/G \}.
\end{align*} 
Using the chi-square identity, we have
\begin{align*}
 {\rm E}  \{(S\phi^2)/( F\si^2)\}
& =   {\rm E}  \{ {n\phi^2/ F} + 2S {\phi^2/(FS)} + 
2S {2\phi } ( - {F}\phi_F/S - {G }\phi_G/S+\phi_S )/F \}
\\
& =
 {\rm E}  \{ (n+2) {\phi^2/F} - 4{\phi }\big(F\phi_F+G\phi_G\big)/F+ 4{S}\phi\phi_S/F \}.
\end{align*}
Similarly,
$$
 {\rm E} \{(S\psi^2)/( G\si^2)\}
=
 {\rm E} \{ (n+2) {\psi^2/ G} - 4{\psi}\big(F\psi_F+G\psi_G\big)/G+ 4{S}\psi\psi_S/G \}.
$$
Thus, one gets
\begin{equation}
I_4 =
 {\rm E} \{ (n+2) {\phi^2/F} - 4{\phi}\big(F\phi_F+G\phi_G\big)/F+ 4{S}\phi\phi_S/F +
 (n+2) {\psi^2/ G} - 4{\psi}\big(F\psi_F+G\psi_G\big)/G+ 4{S}\psi\psi_S/G \}.
\label{eqn:I4}
\end{equation}
Combining (\ref{eqn:I2}), (\ref{eqn:I3}) and (\ref{eqn:I4}) gives the expression in Theorem~\ref{thm:UER}.\hfill $\Box$

\subsection{Proof of Theorem~\ref{thm:HB2}}

It suffices to show that $\phi^{HB2}(F,G,S)$ satisfies conditions (a) and (b) in Proposition~\ref{prp:1}, and that $\psi^{HB2}(F,G,S)$ satisfies conditions (a) and (b) in Proposition~\ref{prp:2}.
For simplicity, let $h(x,y,v)=x^{\al }y^{\be}v^{\ga }\exp\{ - v(x+y+1)/2\}$ for $\al=p(k-1)+a-1$, $\be=p/2+b-1$ and $\ga=(n+pk)/2-c$.
Then, $\phi^{HB2}(F,G,S)$ and $\psi^{HB2}(F,G,S)$ are written as 
\begin{align*}
\phi^{HB2}(F,G,S)=&
{\int_0^F \int_0^G \int_{L  S}^\infty x h(x,y,v) dv dy dx
\Big/
\int_0^F \int_0^G \int_{L  S}^\infty h(x,y,v) dv dy dx},
\\
\psi^{HB2}(F,G,S)=&
{\int_0^F \int_0^G \int_{L  S}^\infty y h(x,y,v) dv dy dx
\Big/
\int_0^F \int_0^G \int_{L  S}^\infty h(x,y,v) dv dy dx}.
\end{align*}
We check conditions (a) and (b) in Proposition~\ref{prp:1} for $\phi^{HB2}(F,G,S)$.
The proof for $\psi^{HB2}(F,G,S)$ is omitted, because it can be shown similarly.

We begin by showing that $\phi^{HB2}(F,G,S)$ is increasing in $F$.
The derivative of $\phi^{HB2}(F,G,S)$ with respect to $F$ is proportional to
\begin{multline*}
\int_0^G \int_{L  S}^\infty F h(F,y,v)dvdy\times \int_0^F\int_0^G\int_{L  S}^\infty  h(x,y,v)dvdydx\\
-\int_0^F\int_0^G\int_{L  S}^\infty x h(x,y,v)dvdydx\times \int_0^G\int_{L  S}^\infty h(F,y,v)dvdy\\
= \int_0^F(F-x) \left\{\int_0^G\int_{L  S}^\infty h(F,y,v)dvdy\times \int_0^G\int_{L  S}^\infty h(x,y,v)dvdy\right\} dx,
\end{multline*}
which is non-negative.
  
We next show that $\phi^{HB2}(F,G,S)$ is increasing in $G$.
The derivative of $\phi^{HB2}(F,G,S)$ with respect to $G$ is proportional to
\begin{multline*}
\int_0^F \int_{L  S}^\infty x h(x,G,v)dvdx\times \int_0^F\int_0^G\int_{L  S}^\infty  h(x,y,v)dvdydx\\
-\int_0^F\int_0^G\int_{L  S}^\infty x h(x,y,v)dvdydx\times \int_0^F\int_{L  S}^\infty h(x,G,v)dvdx\\
= \left\{\int_0^F\int_{L  S}^\infty h(x,G,v)dvdx \right\}^2
\\
\times
\left[  {\rm E}  ^*(X)  {\rm E}  \left\{ {\int_0^G\int_{L  S}^\infty  h(X,y,v)dvdy \over \int_{L  S}^\infty h(X,G,v)dv} \right\} -  {\rm E}  ^*\left\{ X {\int_0^G\int_{L  S}^\infty  h(X,y,v)dvdy \over \int_{L  S}^\infty h(X,G,v)dv}\right\} \right],
\end{multline*}
where $ {\rm E}  ^*$ denotes expectation with respect to the probability 
$$
 {\rm Pr}  (X\in A)=\int_{A\cap[0,F]} \int_{L  S}^\infty h(x,G,v)dvdx\Big/ \int_0^F\int_{L  S}^\infty h(x,G,v)dvdx.
$$
Let $g(X)=\int_0^G\int_{L  S}^\infty  h(X,y,v)dvdy / \int_{L  S}^\infty h(X,G,v)dv$.
If $g(X)$ is decreasing in $X$, it is seen that $g(X)$ and $X$ are monotone in opposite directions, which implies that $ {\rm E}  ^*\{X g(X)\} -  {\rm E}  ^*(X)  {\rm E}  ^*\{g(X)\}\leq 0$ from the covariance inequality.
Thus, we need to show that $g(X)$ is decreasing in $X$.
Note that $g(x)$ is written as
$$
g(x) = \int_0^G\int_{LS}^\infty y^\be v^\ga e^{-v(x+y+1)/2}dvdy\Big/ \int_{LS}^\infty G^\be v^\ga e^{-v(x+G+1)/2}dv.
$$
The derivative $g'(x)$ is proportional to
\begin{align*}
- \int_0^G \int_{LS}^\infty y^\be & v^{\ga+1} e^{-v(x+y+1)/2}dvdy \int_{LS}^\infty v^\ga e^{-v(x+G+1)/2}dv\\
& + \,\int_0^G\int_{LS}^\infty y^\be v^\ga e^{-v(x+y+1)/2}dvdy \int_{LS}^\infty v^{\ga+1} e^{-v(x+G+1)/2}dv\\
\\ &=
\left\{\int_{LS}^\infty v^\ga e^{-v(x+G+1)/2}dv\right\}^2
\left[ -  {\rm E}  ^\dagger \left\{ V {\int_0^G y^\be e^{-V(x+y+1)/2}dy\over e^{-V(x+G+1)/2}}\right\} +  {\rm E}  ^\dagger(V)  {\rm E}  ^\dagger \left\{ {\int_0^G y^\be e^{-V(x+y+1)/2}dy\over e^{-V(x+G+1)/2}}\right\} \right]\\
& =
\left\{\int_{LS}^\infty v^\ga e^{-v(x+G+1)/2}dv\right\}^2
\left[ -  {\rm E}  ^\dagger \left\{ V \int_0^G y^\be e^{(G-y)V}dy\right\} +  {\rm E}  ^\dagger(V)  {\rm E}  ^\dagger \left\{ \int_0^G y^\be e^{(G-y)V}dy \right\} \right],
\end{align*}
where $ {\rm E}  ^\dagger$ denotes the expectation with respect to the probability 
$$
 {\rm Pr}  (V\in A)=\int_{A\cap[LS,\infty)} v^\ga e^{-v(x+G+1)/2}dv\Big/ \int_{LS}^\infty v^\ga e^{-v(x+G+1)/2}dv
$$
 for fixed $x$.
Because $\int_0^G y^\be e^{(G-y)V}dy$ is increasing in $V$ for $y<G$, the covariance inequality implies that
$$
 -  {\rm E}  ^\dagger \left\{ V \int_0^G y^\be e^{(G-y)V}dy\right\} +  {\rm E}  ^\dagger(V)  {\rm E}  ^\dagger \left\{ \int_0^G y^\be e^{(G-y)V}dy\right\}\leq 0,
$$
which means that $g'(x)\leq 0$.
Hence, $\phi^{HB2}(F,G,S)$ is increasing in $G$.
  
We now show that $\phi^{HB2}(F,G,S)$ is decreasing in $S$.
The derivative of $\phi^{HB2}(F,G,S)$ with respect to $S$ is proportional to
\begin{align*}
-L \int_0^F \int_0^G x h(x,y,L S)dydx & \int_0^F\int_0^G\int_{L  S}^\infty  h(x,y,v)dvdydx \\
 +\, L\int_0^F\int_0^G\int_{L  S}^\infty & x h(x,y,v)dvdydx \int_0^F\int_0^G h(x,y,L  S)dydx
\\ & \hspace{-5mm} = L \left\{\int_0^F\int_0^G h(x,y,L  S)dydx\right\}^2 \times \\
& \quad \quad \left[ -  {\rm E}  ^{**}(X)\times  {\rm E}  ^{**} \left\{ { \int_0^G\int_{L  S}^\infty  h(X,y,v)dvdy \over \int_0^G h(X,y,L  S)dy} \right\} +  {\rm E}  ^{**} \left\{  X { \int_0^G\int_{L  S}^\infty  h(X,y,v)dvdy \over \int_0^G h(X,y,L  S)dy}\right\} \right],
\end{align*}
where $ {\rm E}  ^{**}$ is the expectation with respect to the probability
$$
 {\rm Pr}  (X\in A) = \int_{A\cap[0,F]} \int_0^G h(x,y,L  S)dydx\Big/ \int_0^F\int_0^G h(x,y,L  S)dydx.
$$
Note that
$$
{ \int_0^G\int_{L  S}^\infty  h(X,y,v)dvdy} {\Big /} { \int_0^G h(X,y,L  S)dy}
= {\int_0^G\int_{LS}^\infty y^\be v^\ga e^{(LS-v)X - v(y+1)/2}dvdy} {\Big /} { \int_0^G y^\be(LS)^\ga e^{-LS(y+1)/2}dy},
$$
which is decreasing in $X$.
Thus, $X$ and $\int_0^G\int_{L  S}^\infty  h(X,y,v)dvdy /\int_0^G h(X,y,L  S)dy$ are monotone in opposite directions, so that the derivative $\phi^{HB2}(F,G,S)$ with respect to $S$ is negative due to the covariance inequality.
Hence, condition (a) is satisfied for $\phi^{HB2}(F,G,S)$.

Finally, from condition (a), we see that
\begin{align*}
\phi^{HB2}(F,G,S)\leq& \lim_{F\to\infty}\lim_{G\to\infty}\lim_{S\to 0} \phi^{HB2}(F,G,S)
=
{\int_0^\infty\int_0^\infty \int_{0}^\infty x^{\al+1 }y^{\be }v^\ga\exp\{ - y(x+y+1)/2\}dvdy dx
\over
\int_0^\infty \int_{0}^\infty x^{\al}y^{\be }v^\ga\exp\{ - y(x+y+1)/2\}dvdy dx}\\
=& 
{\Ga(\al+2)\Ga(\be+1)\Ga(\ga-\al-\be-2)\over \Ga(\al+1)\Ga(\be+1)\Ga(\ga-\al-\be-1)}
={\al+1 \over \ga-\al-\be-2}={p(k-1)+2a \over n- 2(a+b+c)}.
\end{align*}
Similarly,
\begin{align*}
\psi^{HB2}(F,G,S)\leq& \lim_{F\to\infty}\lim_{G\to\infty}\lim_{S\to 0} \psi^{HB2}(F,G,S)
=
{\int_0^\infty\int_0^\infty \int_{0}^\infty x^{\al }y^{\be+1 }v^\ga\exp\{ - y(x+y+1)/2\}dvdy dx
\over
\int_0^\infty \int_{0}^\infty x^{\al}y^{\be }v^\ga\exp\{ - y(x+y+1)/2\}dvdy dx}\\
=& 
{\Ga(\al+1)\Ga(\be+2)\Ga(\ga-\al-\be-2)\over \Ga(\al+1)\Ga(\be+1)\Ga(\ga-\al-\be-1)}
={\be+1 \over \ga-\al-\be-2}={p+2b \over n- 2(a+b+c)}.
\end{align*}
Thus, condition (b) in Proposition~\ref{prp:1} is satisfied if $a +b +c < n/2$ and if 
$$
\{p(k-1)+2a\}/\{n- 2(a+b+c)\} \leq 2 \{p(k-1)-2 \}( n+2),
$$
which is equivalently rewritten by (\ref{eqn:HB2c}).
Also, condition (b) in Proposition~\ref{prp:2} is satisfied if $a +b +c < n/2$ and if 
$$
(p+2b)/\{n- 2(a+b+c)\} \leq 2(p-2)/(n+2),
$$
which leads to (\ref{eqn:HB2cc}).
Hence, the minimaxity of the hierarchical Bayes estimator is established.

\section*{Acknowledgments}

We thank the Editor-in-Chief, the Associate Editor and the reviewers for valuable comments and helpful suggestions which led to an improved version of this paper. Research of the first author was supported by Grant-in-Aid for JSPS Research Fellows (18J21162). Research of the second author was supported in part by Grant-in-Aid for Scientific Research  (18K11188, 15H01943, 26330036) from Japan Society for the Promotion of Science.


\end{document}